\documentclass[12pt, twoside]{article}

\usepackage{amsfonts}
\usepackage{amssymb}
\usepackage{graphics, epic, eepic}
\usepackage{amsmath, amsthm, latexsym}

\def \c{\mathbb{C}}
\def \z{\mathbb{Z}}
\def \r{\mathbb{R}}

\def \p{\mathbb{P}}

\def \tor{(\c^{*})^d}

\def \e{\mathcal{E}}
\def \o{\mathcal{O}}
\def \Z{\mathcal{Z}}
\def \V{\mathcal{V}}

\def \.{\cdot}

\def \dim{\textup{dim}}

\def \id{\textup{id}}

\def \max{\textup{max}}
\def \Span{\textup{Span}}

\theoremstyle{plain}
\newtheorem{Th}{Theorem}[section]

\newtheorem{Prop}[Th]{Proposition}

\theoremstyle{definition}
\newtheorem{Ex}{Example}[section]

\newtheorem{Rem}{Remark}[section]

\begin{document}
\title{Fixed Points of Torus Action and Cohomology Ring of Toric Varieties}
\author{Kiumars Kaveh\\Department of Mathematics\\University of British Columbia\\} \maketitle

{\footnotesize {\bf Abstract.} Let $X$ be a smooth simplicial
toric variety. Let $Z$ be the set of $T$-fixed points of $X$. We
construct a filtration $F_0 \subset F_1 \subset \cdots$ of $A(Z)$,
the ring of $\c$-valued functions on $Z$, such that $Gr A(Z) \cong
H^*(X, \c)$ as graded algebras. This is the explanation of the
general results of Carrell and Liebermann on the cohomology of
$T$-varities, in the case of toric varieties. We give an explicit
isomorphism between $Gr A(Z)$ and Brion's description of the
polytope algebra.
}


\tableofcontents

\section{Introduction}
In \cite{Carr-Lieb} and \cite{Carr-Lieb1}, Carrell and Lieberman
prove that if $X$ is a smooth projective variety over $\c$ with a
holomorphic vector field $\V$ such that the $Zero(\V)$ is
non-trivial and isolated, then the coordinate ring $A(Z)$ of the
zero scheme $Z$ of $\V$ admits a filtration $F_0 \subset F_1
\subset \cdots$ such that the associated graded $Gr A(Z)$ is
isomorphic to $H^*(X, \c)$ as graded algebra. In this paper, we
give an explicit construction of this filtration in the toric
case. We give an explicit isomorphism between $Gr A(Z)$ and
Brion's description of the polytope algebra (see \cite{Brion}). We
also give direct proofs that the usual relations in the cohomolgy
of a toric variety hold in $Gr A(Z)$.

In the toric case, for the vector field $\V$ one takes the
generating vector field of a $1$-parameter subgroup $\gamma$ in
general position of the torus $T$, so that the fixed point set of
$\gamma$ is the same as the fixed point set of $T$.

This paper is motivated in part by a comment of T. Oda. In
\cite{Oda} p.417, Oda comments about how to explain the results of
Carrell-Lieberman in the toric case: as Khovanskii has shown in
\cite{Askold2}, composition of $\gamma$ and the moment map of the
toric variety $X$ defines a Morse function on $X$ whose critical
points are the fixed points (see Remark \ref{rem-Morse}). Since
the number of critical points of index $i$ is the $i$-th Betti
number, Oda reasonably suggests that the grading on the fixed
point set induced by the Morse index is the grading in
Carrell-Lieberman and hence gives the cohomology algebra. It
happens that this is not necessarily correct. One can see this in
the example of $\c \p ^2$ (see Remark
\ref{rem-Oda-counter-example}).


\noindent{\bf Acknowledgement:} I would like to express my
gratitude to Prof. James Carrell for suggesting the problem to me
and helpful discussions. I would also like to thank
my friend, Vladlen Timorin for helpful discussions.

\section{Preliminaries on Cohomology of Varieties with $G_m$ Action}
\label{sec-cohomology-torus-action} In this section we briefly
review the general theorems due to Carrell and Lieberman (see
\cite{Carrell} and \cite{Carr-Lieb1}) on the cohomology of
varieties with a $G_m$ action.

Let $X$ be a smooth projective variety over $\c$.

\begin{Th}[\cite{Carrell}, Theorem 5.4] \label{theorem-filtration}
Suppose $X$ admits a holomorphic vector field $\V$ with $Zero(\V)$
isolated but non-trivial. Then the coordinate ring $A(Z)$ of the
zero scheme $Z$ of $\V$ admits an increasing filtration $F_\bullet
= F_\bullet A(Z)$ such that
\begin{itemize}
\item[(i)] $F_iF_j \subset F_{i+j}$; and
\item[(ii)] $H^*(X, \c) = \bigoplus_{i\geq0}
H^{2i}(X, \c) \cong \bigoplus_{i\geq 0}Gr_{2i}(A(Z))$,
\end{itemize}
where the displayed summands are isomorphic over $\c$. Here
$$Gr_{2i}(A(Z)) := F_iA(Z) / F_{i-1}A(Z).$$
\end{Th}

Let $E \rightarrow X$ be a holomorphic vector bundle and $\e$ its
sheaf of holomorphic sections. One says that $E$ is
$\V$-equivariant if the derivation $\V$ of $\o _x$ lifts to $\e$.
That is, there exists a $\c$-linear sheaf homomorphism
$\tilde{\V}: \e \rightarrow \e$ such that if $\sigma \in \e_x$ and
$f \in \o_{X,x}$ then
$$\tilde{\V}(f\sigma) = \V(f)\sigma + f\tilde{\V}(\sigma).$$
We then have:
\begin{Th}[\cite{Carrell}, Theorem 5.5] \label{theorem-chern-in-Gr}
If $p$ is a polynomial of degree $l$, then $p(\tilde{\V}_{|Z}) \in
F_lA(Z)$, and in the associated graded, that is, in $Gr_{2l}A(Z)$,
$p(\tilde{\V}_{|Z})$ corresponds to $p(c(\e)) \in H^l(X, \Omega^l)
= H^{2l}(X, \c)$, where $c(\e)$ denotes the Atiyah-Chern class of
$\e$.
\end{Th}

\section{Preliminaries on the Cohomology of Toric Varieties}
\label{sec-cohomology-toric-var} Let $T$ be the algebraic torus
$\tor$. As usual, $N$ denotes the lattice of $1$-parameter
subgroups of $T$, $N_\r$ the real vector space $N \otimes_\z \r$,
$M$ the dual lattice of $N$ which is the lattice of characters of
$T$ and, $M_\r$ the real vector space $M \otimes_\z \r$. A vector
$n = (n_1,\ldots,n_d) \in \z^d \cong N$ corresponds to the
$1$-parameter subgroup $t^n = (t^{n_1},\ldots, t^{n_d})$ .
Similarly, a covector $m = (m_1, \ldots, m_d) \in (\z^d)^* \cong
M$ corresponds to the character $x^m = x_1^{m_1}\cdots x_d^{m_d}$.
We use $\langle \phantom{a},\phantom{a}\rangle:N \times M
\rightarrow \z$ for the natural pairing between $N$ and $M$.

Let $X$ be a $d$-dimensional smooth projective simplicial toric
variety. Let $\Sigma \subset N_\r$ be the simplicial fan
corresponding to $X$. We denote by $\Sigma(i)$ the set of all
$i$-dimensional cones in $\Sigma$. For each $\rho \in \Sigma(1)$,
let $\xi_\rho$ be the primitive vector along $\rho$, i.e. the
smallest integral vector on $\rho$.

There is a 1-1 correspondence between the
orbits of dimension $i$ in $X$ and the cones in $\Sigma(d-i)$.
The fixed points of $T$ correspond to
the cones in $\Sigma(d)$.
In a smooth toric variety all the orbit closures are
smooth, the cohomology class dual to the closure of the orbit
corresponding to $\rho \in \Sigma(1)$ is denoted by $D_\rho \in H^2(X, \c)$.
It is well-known that the cohomology algebra of a toric variety is
generated by the classes $D_\rho$. More precisely, we have:

\begin{Th}[see \cite{Fulton}, p.106]
\label{theorem-cohomology-toric-var}
Let $X$ be a smooth projective toric variety. Then $H^*(X, \c) =
\z[D_\rho, \rho \in \Sigma(1)] / I$, where $I$ is the ideal
generated by all
\begin{itemize}
\item[(i)] $D_{\rho_1}\. \ldots \. D_{\rho_k},
\quad \forall \rho_1, \ldots, \rho_k$ not in a cone of $\Sigma$; and
\item[(ii)] $\sum_{\rho \in \Sigma(1)} \langle \xi_\rho, u \rangle
D_\rho,  \quad \forall u \in M$.
\end{itemize}
\end{Th}

Now, let $\Delta \subset M_\r$ be a simple rational polytope
normal to the fan $\Sigma$. The polytope $\Delta$ defines a
diagonal representation $\pi: T \rightarrow GL(V)$ where
$\dim_\c(V) = $ the number of lattice points in $\Delta$. If the
mutual differences of the lattice points in $\Delta$ generate $M$
then we get an embedding of $X$ in $\p (V)$ as the closure of the
orbit of $(1:\ldots:1)$. In the rest of the paper, we assume that
the above condition holds for $\Delta$.

The set of faces of dimension $i$ in $\Delta$ is denoted by
$\Delta(i)$. There is a 1-1 correspondence between the the faces
in $\Delta(i)$ and the cones in $\Sigma(d-i)$ which in turn
correspond to the orbits of dimension $i$ in $X$. Hence the fixed
points of $T$ on $X$ correspond to the vertices of $\Delta$.

The support function $l_\Delta: N_\r \rightarrow \r$ is
defined by: $l_\Delta(\xi) = \max_{x\in\Delta} \langle \xi, x
\rangle$.

Let $L_\Delta$ be the line bundle on $X$ obtained by restricting
the universal subbundle on $\p (V)$ to $X$. We will need the
following classical theorem which tells us how the first Chern
class $c_1(L_\Delta)$ is represented as a linear combination of
the classes $D_\rho$.

\begin{Th}
With notation as above we have
$$c_1(L_\Delta) = \sum_{\rho \in
\Sigma(1)} l_\Delta(\xi_\rho) D_\rho.$$
\end{Th}


\section{Main Theorem}
\label{sec-main-theorem} As before, let $X$ be a smooth projective
simplicial toric variety with fan $\Sigma$ and a polytope $\Delta$
normal to the fan which gives rise to a representation $\pi: T
\rightarrow GL(V)$ and a $T$-equivariant embedding of $X$ in $\p
(V)$, for a vector space $V$ over $\c$. Let $\gamma \in N$ be a
$1$-parameter subgroup of $T$. We can choose $\gamma$ so that the
set of fixed points of $\gamma$ is the same as the set of fixed
points of $T$. We denote the set of fixed points by Z.

In this section, we construct a filtration $F_0 \subset F_1
\subset \cdots$ for $A(Z)$ such that $H^*(X, \c) \cong Gr A(Z)$.

\noindent{\bf Notation:} In the following, $z$ denotes a fixed
point, $\sigma$ the corresponding $d$-dimensional cone in $\Sigma$
and $v$ the corresponding vertex in $\Delta$. A $1$-dimensional
cone in $\Sigma$ is denoted by $\rho$ and the corresponding facet
of $\Delta$ by $F$.

From Theorem \ref{theorem-filtration} applied to the generating
vector field of $\gamma$, there exists a filtration $F_0 \subset
F_1 \subset \cdots$ of $A(Z)$, the ring of $\c$ valued functions
on $Z$, so that $H^*(X, \c) \cong \bigoplus_{i=0}^\infty F_{i+1} /
F_{i}$, as graded algebras. In particular, we have $H^2(X, \c)
\cong F_1 / F_0$. The subspace $F_0 A(Z)$ is just the set of
constant functions. To determine the image of $H^2(X, \c)$ in $Gr
A(Z)$ we need to determine $F_1$. We start by finding the
representatives in $F_1$ for the Chern classes of the line
bundles.

The $1$-parameter subgroup $\gamma: \c^* \rightarrow T$ acts on
$V$ via $\pi$ and hence the action of $\gamma$ on $X$ lifts to an
action of $\gamma$ on the line bundle $L_\Delta$. Thus the
generating vector field of $\gamma$ has a lift to $L_\Delta$. If
we view $L_\Delta$ as $\{(x, l) \in X \times V \mid x = [l] \}$
then the action of $\gamma$ on $L_\Delta$ is given by:
$$\gamma(t) \. (x,l) = ( \pi(t^\gamma)x, \pi(t^\gamma)l).$$
Now, from Theorem \ref{theorem-chern-in-Gr} we have:
\begin{Prop} \label{prop-f_Delta}
Under the isomorphism $F_1 / F_0 \cong H^2(X, \c)$, the first
Chern class $c_1(L_\Delta)$ is represented by the function
$f_\Delta$ defined by:
$$ f_\Delta(z) = \langle \gamma, v \rangle, \quad \forall z \in
Z,$$ where $v$ is the vertex of $\Delta$ corresponding to the
fixed point $z$.
\end{Prop}
\begin{proof}
In Theorem \ref{theorem-chern-in-Gr}, take $E$ to be $L_\Delta$
and $p$ be the identity polynomial. The derivation $\tilde{\V}$ is
just the derivation given by the $G_m$-action of $\gamma$ on
$L_\Delta$. Let $z$ be a fixed point and $(z,l) \in (L_\Delta)_z$
a point in the fiber of $z$. We have:
\begin{eqnarray*}
\gamma(t)\. (z, v) &=& (z, \pi(t^\gamma)l), \cr &=& (z, \langle
\gamma , v \rangle l). \cr
\end{eqnarray*}
and hence $f_\Delta(z) = \langle \gamma , v \rangle$.
\end{proof}

\begin{figure}[ht]
\begin{center}
\setlength{\unitlength}{0.00083333in}
\begingroup\makeatletter\ifx\SetFigFont\undefined%
\gdef\SetFigFont#1#2#3#4#5{%
  \reset@font\fontsize{#1}{#2pt}%
  \fontfamily{#3}\fontseries{#4}\fontshape{#5}%
  \selectfont}%
\fi\endgroup%
{\renewcommand{\dashlinestretch}{30}
\begin{picture}(3184,3602)(0,-10)
\path(1670,3106)(178,2609)(12,950)
	(2057,12)(2885,1779)(1670,3106)
\dashline{90.000}(2885,1779)(3162,2388)
\path(1504,1503)(2223,2222)
\blacken\thicklines
\path(2129.227,2065.712)(2223.000,2222.000)(2066.712,2128.227)(2129.227,2065.712)
\put(1615,3271){\makebox(0,0)[lb]{\smash{{{\SetFigFont{9}{10.8}{\familydefault}{\mddefault}{\updefault}$v$}}}}}
\put(2333,3492){\makebox(0,0)[lb]{\smash{{{\SetFigFont{9}{10.8}{\familydefault}{\mddefault}{\updefault}$v'$}}}}}
\put(2167,2719){\makebox(0,0)[lb]{\smash{{{\SetFigFont{9}{10.8}{\familydefault}{\mddefault}{\updefault}$F$}}}}}
\put(2831,2995){\makebox(0,0)[lb]{\smash{{{\SetFigFont{9}{10.8}{\familydefault}{\mddefault}{\updefault}$F'$}}}}}
\put(1891,1668){\makebox(0,0)[lb]{\smash{{{\SetFigFont{9}{10.8}{\familydefault}{\mddefault}{\updefault}$\xi_\rho$}}}}}
\put(2555,563){\makebox(0,0)[lb]{\smash{{{\SetFigFont{9}{10.8}{\familydefault}{\mddefault}{\updefault}$\Delta$}}}}}
\path(2278,3327)(3162,2388)
\path(1670,3106)(2885,1779)
\thinlines
\dashline{90.000}(1670,3106)(2278,3327)
\end{picture}
}
\caption{Moving facet $F$} ~\label{Delta}
\end{center}
\end{figure}
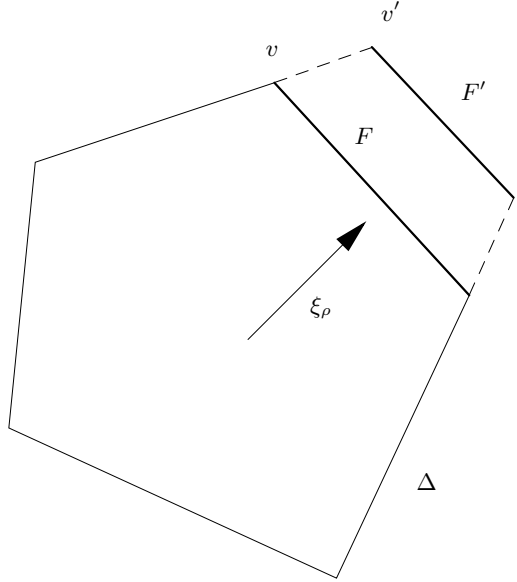

Next, we wish to determine the images of the classes $D_\rho, \rho
\in \Sigma(1)$, in $F_1 /F_0$. Fix a $1$-dimensional cone $\rho$
in $\Sigma(1)$. Let $F$ be the facet of $\Delta$ orthogonal to
$\rho$. We move the facet $F$ of $\Delta$ parallely to obtain a
new polytope $\Delta'$ (Figure \ref{Delta}). The polytope
$\Delta'$ is still normal to the fan $\Sigma$. Let $F'$ denote the
facet of $\Delta'$ obtained by moving $F$. The maximum of the
function $\langle \xi_\rho , \. \rangle$ on $\Delta$ and $\Delta'$
is obtained on the facets $F$ and $F'$ respectively. For support
functions of these polytopes we can write:
$$ l_\Delta(\xi_\rho) = \langle \xi_\rho, \text{ some point in }F \rangle,$$
$$ l_{\Delta'}(\xi_\rho) = \langle \xi_\rho, \text{ some point in }F' \rangle$$
$$ l_\Delta(\xi_{\rho'}) = l_{\Delta'}(\xi_{\rho'}), \quad \forall \rho' \neq \rho.$$
We also have:
$$ c_1(L_\Delta) = l_\Delta(\xi_\rho) D_\rho +
\sum_{\rho' \in \Sigma(1), \rho' \neq \rho} l_\Delta(\xi_{\rho'}) D_{\rho'},$$
$$ c_1(L_{\Delta'}) = l_{\Delta'}(\xi_\rho) D_\rho +
\sum_{\rho' \in \Sigma(1), \rho' \neq \rho}
l_{\Delta'}(\xi_{\rho'}) D_{\rho'}.$$ Hence
$$c_1(L_\Delta) - c_1(L_{\Delta'}) = (l_\Delta(\xi_\rho) -
l_{\Delta'}(\xi_\rho)) D_\rho.$$ So
$$ D_\rho = \frac{c_1(L_\Delta) - c_1(L_{\Delta'})}{l_\Delta(\xi_\rho) -
l_{\Delta'}(\xi_\rho)}.$$


Now, let $z$ be a fixed point, $\sigma$ the corresponding
$d$-dimensional cone, and $v$ and $v'$ the corresponding vertices
in $\Delta$ and $\Delta'$ respectively. From Proposition
\ref{prop-f_Delta}, $D_\rho$ corresponds to the function $f_\rho
\in F_1 A(Z)$ given by:
\begin{eqnarray*}
f_\rho(z) &=& \frac{f_\Delta(z) -
f_{\Delta'}(z)}{l_\Delta(\xi_\rho) - l_{\Delta'}(\xi_\rho)}, \cr
&=& \frac{\langle \gamma , v-v' \rangle}{l_\Delta(\xi_\rho) -
l_{\Delta'}(\xi_\rho)}. \cr
\end{eqnarray*}

If $v \notin F$ then $v=v'$ and hence $f_\rho(z) = 0$. If $v \in
F$ then $l_\Delta(\xi_\rho) = \langle \xi_\rho , v \rangle$ and
$l_{\Delta'}(\xi_\rho) = \langle \xi_\rho , v' \rangle$. We obtain
that:
$$ f_\rho(z) = \begin{cases} \frac{\langle \gamma ,
v - v' \rangle}{\langle \xi_\rho, v-v' \rangle} \quad
\text{if } v \in F \\ 0 \quad \text{if } v \notin F \\
\end{cases}$$

Since $\Delta$ is a simple polytope, there are $d$ edges at the
vertex $v$. If $v \in F$, then there is only one edge $e$ at $v$
which does not belong to $F$. The vector $v-v'$, in fact, is along
this edge. Note that the above formula for $f_\rho(z)$ does not
depend on the length of the vector $v-v'$ (i.e. how much we move
the facet $F$ to obtain the new polytope $\Delta'$). Let
$u_{\sigma, \rho}$ be the vector along the edge $e$ normalized
such that $\langle u_{\sigma, \rho}, \xi_\rho \rangle = 1$. Then
we have:
\begin{Prop} \label{prop-f-rho}
With notation as above, the cohomology class $D_\rho$ is
represented by the function $f_\rho$ in $F_1 A(Z)$ defined by
$$ f_\rho(z) = \begin{cases} \langle \gamma ,
u_{\sigma, \rho} \rangle \quad
\text{if } v \in F \\ 0 \quad \text{if } v \notin F \\
\end{cases}$$
\end{Prop}

Since $H^2(X, \c)$ is generated by the classes
$D_\rho, \rho \in \Sigma(1)$ and $H^*(X, \c)$ is generated in degree $2$,
from Theorem \ref{theorem-chern-in-Gr} we obtain:
\begin{Th}
$F_1 A(Z) / F_0 A(Z) = \Span_\c \langle f_\rho, \rho \in \Sigma(1) \rangle$.
Moreover, $F_i A(Z) = $ all polynomials of degree $\leq i$ in
the $f_\rho$.
\end{Th}

One can prove directly that the functions $f_\rho, \rho \in
\Sigma(1)$, satisfy the relations in the statement of Theorem
\ref{theorem-cohomology-toric-var}. More precisely:
\begin{Th} The functions $f_\rho, \rho \in \Sigma(1)$, satisfy the
following relations:
\begin{itemize}
\item[(i)] $f_{\rho_1}\. \ldots \. f_{\rho_k} = 0,
\quad \forall \rho_1, \ldots, \rho_k$ not in a cone of $\Sigma$; and
\item[(ii)] $\sum_{\rho \in \Sigma(1)} \langle \xi_\rho, u \rangle
f_\rho =$ some constant function on $Z$, $\quad \forall u \in M$.
\end{itemize}
\end{Th}
\begin{proof}
(i) is easy because every $f_\rho$ is non-zero only at $z$ such
that the corresponding vertex lies in the facet $F_\rho$
corresponding to $\rho$. Now, if $\rho_1, \ldots, \rho_k$ are not
in a cone of $\Sigma$, it means that the the intersection of the
corresponding facets $F_{\rho_i}$ is empty, i.e. the product of
the $f_{\rho_i}$ is zero.

For (ii), let $z$ be a fixed point and, $\sigma$ and $v$ the
corresponding $d$-dimensional cone and vertex respectively. Let
$A$ be the $d \times d$ matrix whose rows are vectors $\xi_\rho$
and let $B$ be the $d \times d$ matrix whose columns are vectors
$u_{\sigma, \rho}$, where $\rho$ is an edge of $\sigma$. Since the
cone at the vertex $v$, which is generated by the vectors
$u_{\sigma, \rho}$, is dual to the cone $\sigma$, we get $AB =
\id$. Now, we have
\begin{eqnarray*}
\sum_{\rho \in \Sigma(1)} \langle \xi_\rho, u \rangle f_\rho &=&
\sum_{\rho \text{ an edge of } \sigma} \langle \xi_\rho, u \rangle
f_\rho \cr &=&  \sum_{\rho \text{ an edge of } \sigma} \langle
\xi_\rho, u \rangle \langle \gamma, u_{v,F_{\rho}} \rangle \cr &=&
A \. u \. \gamma \. B,
\end{eqnarray*}
where $\.$ means product of matrices and $\gamma$ is regarded as a
row vector and $u$ is regarded as a column vector. But the result
of the above is simply $\langle \gamma, u\rangle$, since $AB =
\id$. So we proved that the expression (ii) is independent of $z$
and hence is a constant function on $Z$.
\end{proof}

One can introduce a finite affine set $\Z$ isomorphic to $Z$ such
that the natural grading on the coordinate ring $A(\Z)$ coincides
with the above filtration $F_{\bullet}$ given by the $f_\rho$.
Define the function $\Theta: Z \rightarrow \r^{\Sigma(1)} \subset
\c^{\Sigma(1)}$ by
$$\Theta(z)_\rho = f_\rho(z),$$
and let $\Z = \Theta(Z)$.
\begin{Prop}
$Gr A(\Z) \cong H^*(X, \c)$, as graded algebras. The grading on
$A(\Z)$ is induced from the usual grading of the polynomial
algebra.
\end{Prop}
\begin{proof}
Immediate.
\end{proof}

\begin{Rem} \label{rem-Morse}
Let $\mu: X \rightarrow M_\r$ be the moment map of the toric
variety and, as before, $\gamma \in N$ a $1$-parameter subgroup in
general position. In \cite{Askold2} Khovanskii shows that the
composition of $\gamma$ and $\mu$ defines a Morse function on $X$
whose critical points are the fixed points of $X$. The Morse index
of a fixed point corresponding to a vertex $v$ is twice the number
of edges at $v$ on which the linear function $\gamma$ is
decreasing. Back to the definition of the functions $f_\rho$
(Proposition \ref{prop-f-rho}), the linear function $\gamma$ is
decreasing on the edge $e$ at $v$ if and only if $f_\rho(z) < 0$.
That is, the Morse index of a fixed point $z$ is equal to twice
the number of negative coordinates of the point $\Theta(z) \in
\r^{\Sigma(1)}$. Since the number of critical points of index $2i$
is the $2i$-th Betti number of $X$, we conclude the non-trivial
relation that: the number of points in $\Z$ exactly $i$ of their
coordinates are negative is equal to $\dim Gr_i A(\Z)$.
\end{Rem}


\section{Relation with the Polytope Algebra}
To each simplicial polytope $\Delta$, one can associate an
algebra, called the {\it polytope algebra} of $\Delta$ (see
\cite{Kh-P}, and for a more detailed explanation \cite{Timorin}).
The direct limit of these algebras for all $\Delta$ is the
McMullen's polytope algebra. McMullen's polytope algebra plays an
important role in the study of finitely additive measures on the
convex polytopes. For an integral polytope $\Delta$, its polytope
algebra coincides with the cohomology algebra of the corresponding
toric variety $X$.

In \cite{Brion}, Brion gives a description of the polytope algebra
of a polytope as a quotient of the algebra of continuous piecewise
polynomial functions: let $\Sigma \subset N_\r$ be the fan of the
polytope $\Delta$. Let $R$ be the algebra of all continuous
functions on $N_\r$ which restricted to each cone of $\Sigma$ are
given by a polynomial. Let $I$ be the ideal of $R$ generated by
all the linear functions on $N_\r$, then the polytope algebra of
$\Delta$ is isomorphic to $R / I$.

There is a good set of generators for $R$ parameterized by the set
of $1$-dimensional cones $\Sigma(1)$. For each $\rho \in
\Sigma(1)$, define $g_\rho: N_\r \to \r$ as a piecewise linear
function, supported on the cones containing $\rho$, as follows:
\begin{itemize}
\item[(i)] $g_\rho =0$ on any cone not containing $\rho$; and
\item[(ii)] for a $d$-dimensional cone $\sigma$ containing $\rho$,
the function $g_\rho$ restricted to $\sigma$ is the unique linear
function defined by $g_\rho(x) = 0$ for $x \in \rho' \neq \rho,
\rho' \in \Sigma(1)$ and $g_\rho (\xi_\rho) = 1$.
\end{itemize}

One can show that the $g_\rho$ are a set of generators for $R$.
Moreover, by sending $g_\rho$ to $D_\rho$, we get an isomorphism
between $R/I$ and $H^*(X, \c)$, in particular the $g_\rho$ satisfy
the relations in Theorem \ref{theorem-cohomology-toric-var}.

In the next theorem, we show how this description of the
cohomology is related to the $Gr A(Z)$ description:

Let $\gamma$ be a $1$-parameter subgroup in general position. Let
$p \in R$ be a continuous piecewise polynomial of degree $n$.
Define $$\Phi(p) = \frac{\partial^np}{\partial\gamma^n},$$ where
$\partial^n/\partial\gamma^n$ means $n$ times differentiation in
the direction of the vector $\gamma$. Then $\Phi(p)$ is a constant
function on each $d$-dimensional cone and hence can be viewed as a
function in $A(Z)$. We have:

\begin{Th} \label{thm-relation-with-Brion-desc}
\begin{itemize}
\item[(i)] $\Phi(g_\rho) = f_\rho$; and
\item[(ii)] $\Phi$ induces an isomorphism between $R/I$ and $Gr A(Z)$.
\end{itemize}
\end{Th}
\begin{proof}(i)
Let $\sigma$ be a $d$-dimensional cone containing $\rho$. Since
$\sigma$ is simplicial the set $\{ \xi_{\rho'} \mid \rho' \subset
\sigma \}$ form a basis for $N_\r$. Consider the linear function
$l$ defined by $l(\xi_\rho) = 1$ and $l(\xi_{\rho'}) = 0, \rho'
\subset \sigma$ and $\rho' \neq \rho)$. Let $A$ be the $d \times
d$ matrix whose rows are vectors $\xi_\rho$ and $B$ be the $d
\times d$ matrix whose columns are vectors $u_{\sigma, \rho}$,
where $\rho$ is an edge of $\sigma$. Let $v$ be the vertex of
$\Delta$ corresponding to $\sigma$. The cone at $v$ is dual to
$\sigma$ and hence we have $AB = \id$. View $\gamma$ as a row
vector. The $\gamma$ in the basis $\xi_{\rho'},\rho' \subset
\sigma$ is $\gamma A^{-1} = \gamma B$. Thus, one sees that the
derivative of $l$ along $\gamma$ is equal to the $\rho$-th
component of $\gamma B$. But this is the same as $f_\rho(z)$.

For (ii), note that $\Phi$ is clearly an additive homomorphism.
Let $p$ and $q$ $\in R$ be of degrees $n$ and $m$ respectively.
Since $n+1$-th derivative of $p$ and $m+1$-th derivative of $q$
along $\gamma$ are zero, one can see that $\Phi(pq) =
\Phi(p)\Phi(q)$. Let $p$ be a linear function on $N_\r$, then the
first derivative of $p$ along $\gamma$ is a constant function
hence $\Phi(p)$ is zero in $Gr A(Z)$, i.e. $\Phi$ is well-defined
on $R/I$. Since the $n$-th graded piece of $R$ (respectively
$A(Z)$) is the set of polynomials of degree $n$ in the $g_\rho$
(respectively $f_\rho$), from (i) we see that $\Phi: R/I \to Gr
A(Z)$, is onto. Since $\dim(R/I)_i = \dim H^{2i}(X, \c) = \dim
F_{i+1}/F_i$, it follows that $\Phi$ is 1-1 as well. This finishes
the proof.
\end{proof}


\section{Examles}
\label{sec-examples} In this section we consider two examples in
dimension $2$, namely, $\c P^2$ and $\c P^1 \times \c P^1$. For
each example, we compute the functions $f_\rho$ and the finite
affine set $\Z$.
\begin{Ex}[$X = \c P^2$]
Fan $\Sigma$ of $\c P^2$ is shown in Figure \ref{fan1}. There are
$3$ one dimensional cones denoted by $\rho_1, \rho_2$ and
$\rho_3$. The primitive vectors along the $\rho_i$ are $\xi_1 =
(1,0), \xi_2 = (0,1)$ and $\xi_3 = (-1,-1)$. The vertices of a
normal polytope to the fan are $v_1 = (1,1), v_2 = (-2,1)$ and
$v_3 = (1, -2)$ (see Figure \ref{polytope1}). They correspond to
the three fixed points $z_1, z_2$ and $z_3$. At each vertex, there
are two vectors along the edges. For $v_1$, we take $\{(1,0),
(0,1)\}$, for $v_2$ we take $\{(-1,0), (-1,1)\}$ and finally, for
$v_3$ we take $\{(0,-1), (1,-1)\}$.

\input epsf
\begin{figure}[ht]
\centerline{\epsfysize=4cm \epsffile{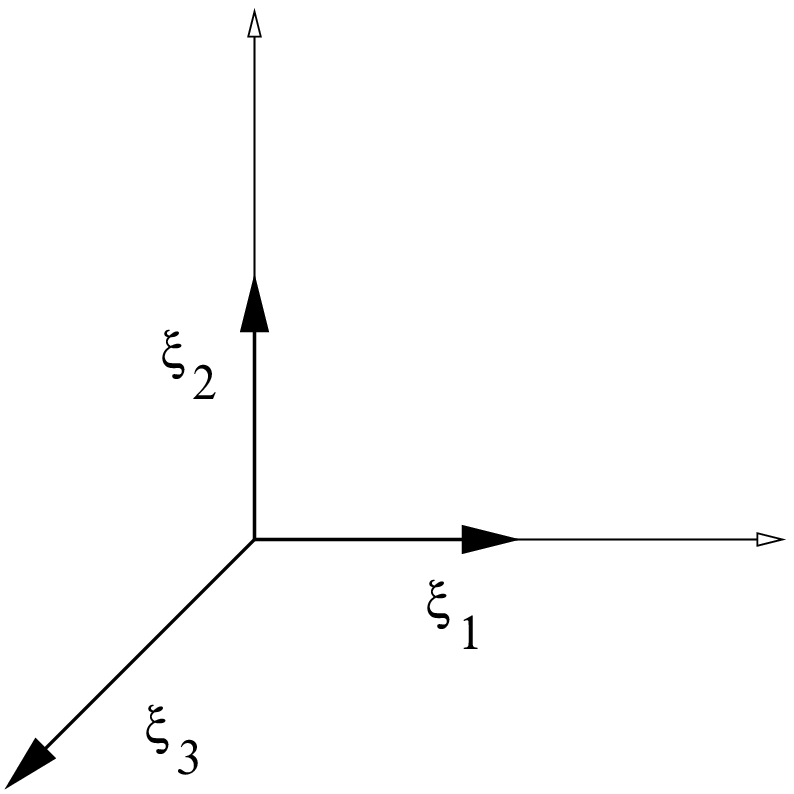}}
\caption{Fan of $\c P^2$} \label{fan1} 
\centerline{\epsfysize=5cm \epsffile{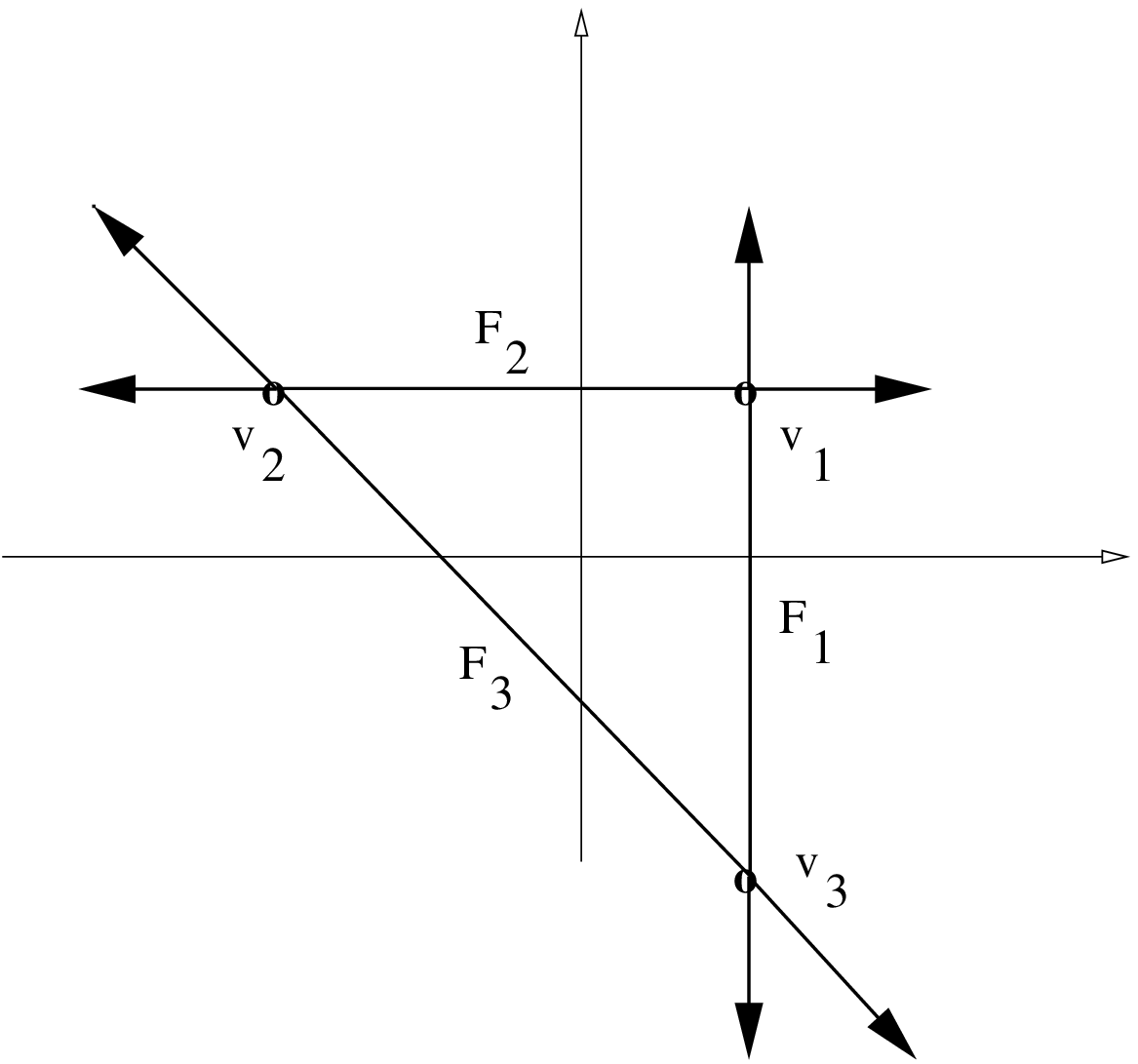}} 
\caption{Polytope normal to the fan of
$\c P^2$ and the vectors $u_{\sigma, \rho}$.} ~\label{polytope1}
\end{figure}

Let $\gamma = (\gamma_1, \gamma_2)$ be a $1$-parameter subgroup.
From the definition of the functions $f_\rho$ (Proposition
\ref{prop-f-rho}), we get the following table for their values:
$$
\begin{array}{|c|c|c|c|}
\hline & & & \\
& z_1 & z_2 & z_3 \\ \hline & & & \\
f_1 & \gamma_2 & \gamma_2 - \gamma_1 & 0 \\ \hline & & & \\
f_2 & \gamma_1 & 0 & \gamma_1 - \gamma_2 \\ \hline & & & \\
f_3 & 0 & -\gamma_1 & -\gamma_2 \\ \hline
\end{array}
$$
and hence, $\Z = \{ (\gamma_2, \gamma_1, 0), (\gamma_2 - \gamma_1,
0, -\gamma_1), (0, \gamma_1 - \gamma_2, -\gamma_2) \} \subset
\r^3$. Note that the points in $\Z$ lie on the same line parallel
to $(1,1,1)$. One can see that $Gr_i A(\Z) \cong \c, 0\leq i\leq
2$ and $Gr_i A(\Z) = \{0\}, i>2$. If $x$ is a non-zero element of
$Gr_1 A(\Z)$ then, $H^*(\c P^2, \c) \cong Gr A(\Z) \cong \c[x] /
\langle x^3 \rangle$.

The above calculation can be carried out in general for $\c P^n$.
One can show that all the points in the set $\Z$ lie on the same
line parallel to $(1,\ldots, 1)$, and hence $Gr_i \cong \c$ for
$0\leq i \leq n$ and $Gr_i \cong 0$ for $i>n$ and thus $H^*(\c
P^n, \c) \cong Gr A(\Z) \cong \c[x] / \langle x^{n+1} \rangle$. In
fact, any set of $n$ points lying on the same line can give the
cohomology of $\c P^n$.

\begin{Rem} \label{rem-Oda-counter-example}
Consider the polytope $\Delta$ for $\c P^2$. For a $\gamma$ in
general position, there is one vertex of index $4$, one vertex of
index $2$ and one vertex of index $0$. Without loss of generality,
assume that the indices of $v_1, v_2$ and $v_3$ are $0, 2$ and $4$
respectively. Now, if the grading on $A(Z)$ is induced by the
Morse index, the subspace of elements of degree $\leq 1$ is
generated by the functions supported on $v_2$, the only fixed
point of index $2$, and the constant functions. Hence, any
function of degree $\leq 1$ should have the same value on $v_1$
and $v_3$. But, for example, $f_1$ does not have this property
while it is of degree $\leq 1$. This shows that the grading by the
Morse index does not coincide with the filtration generated by the
$f_\rho$.
\end{Rem}
\end{Ex}
\begin{Ex}[$X = \c P^1 \times \c P^1$]
Fan $\Sigma$ of $\c P^1 \times \c P^1$ is shown in Figure
\ref{fan2}. There are $4$ one dimensional cones denoted by
$\rho_1, \rho_2, \rho_3$ and $\rho_4$. The primitive vectors along
the $\rho_i$ are $\xi_1 = (1,0), \xi_2 = (0,1), \xi_3 = (-1,0)$
and $\xi_4 = (0,-1)$. The vertices of a normal polytope to the fan
are $v_1 = (1,1), v_2 = (-1,1), v_3 = (-1, -1)$ and $v_4 = (1,
-1)$ (see Figure \ref{polytope2}). They correspond to the four
fixed points $z_1, z_2, z_3$ and $z_4$. At each vertex, there are
two vectors along the edges. For $v_1$, we take $\{(1,0),
(0,1)\}$, for $v_2$ we take $\{(-1,0), (0,1)\}$, for $v_3$ we take
$\{(-1,0),(0,-1)\}$ and finally, for $v_4$ we take $\{(1,0),
(0,-1)\}$.

\input epsf
\begin{figure}[ht]
\centerline{\epsfysize=4cm \epsffile{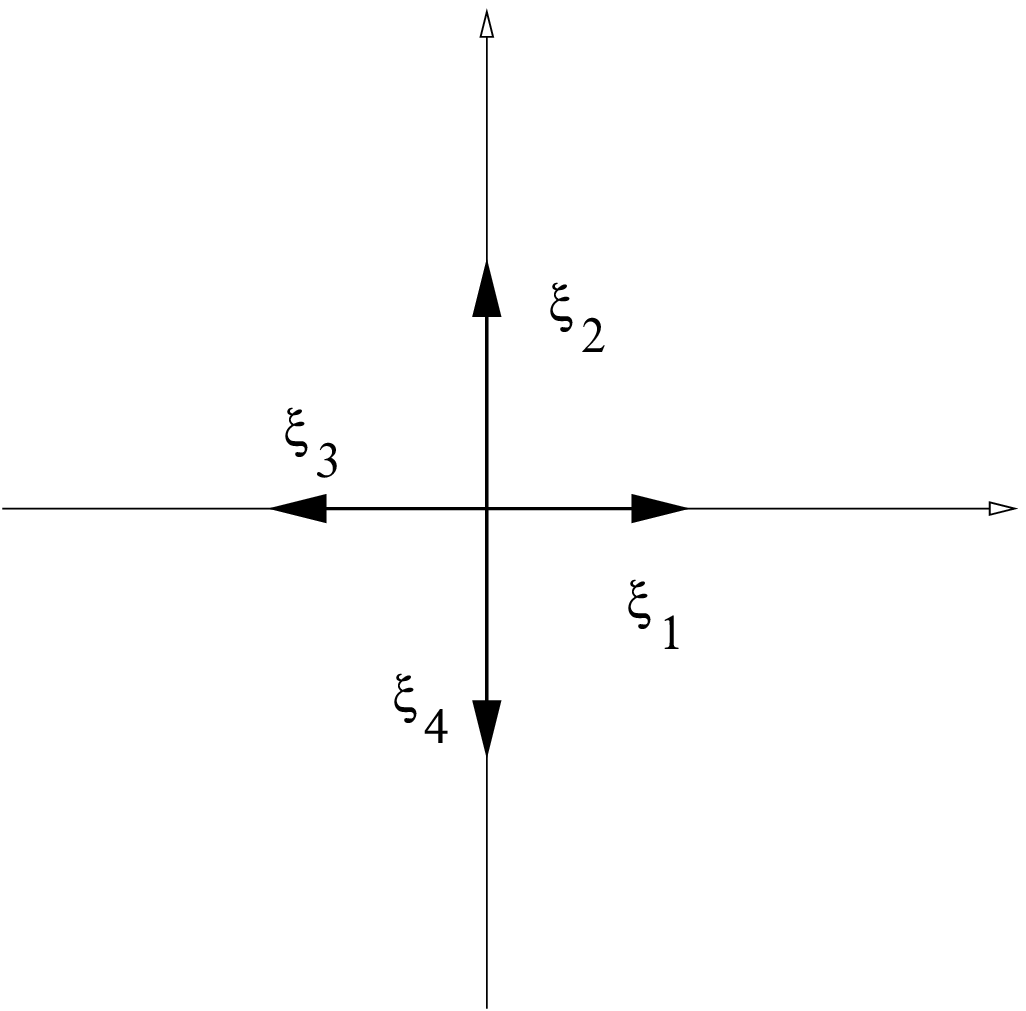}} 
\caption{Fan of $\c P^1 \times \c P^1$} ~\label{fan2} 
\centerline{\epsfysize=5cm \epsffile{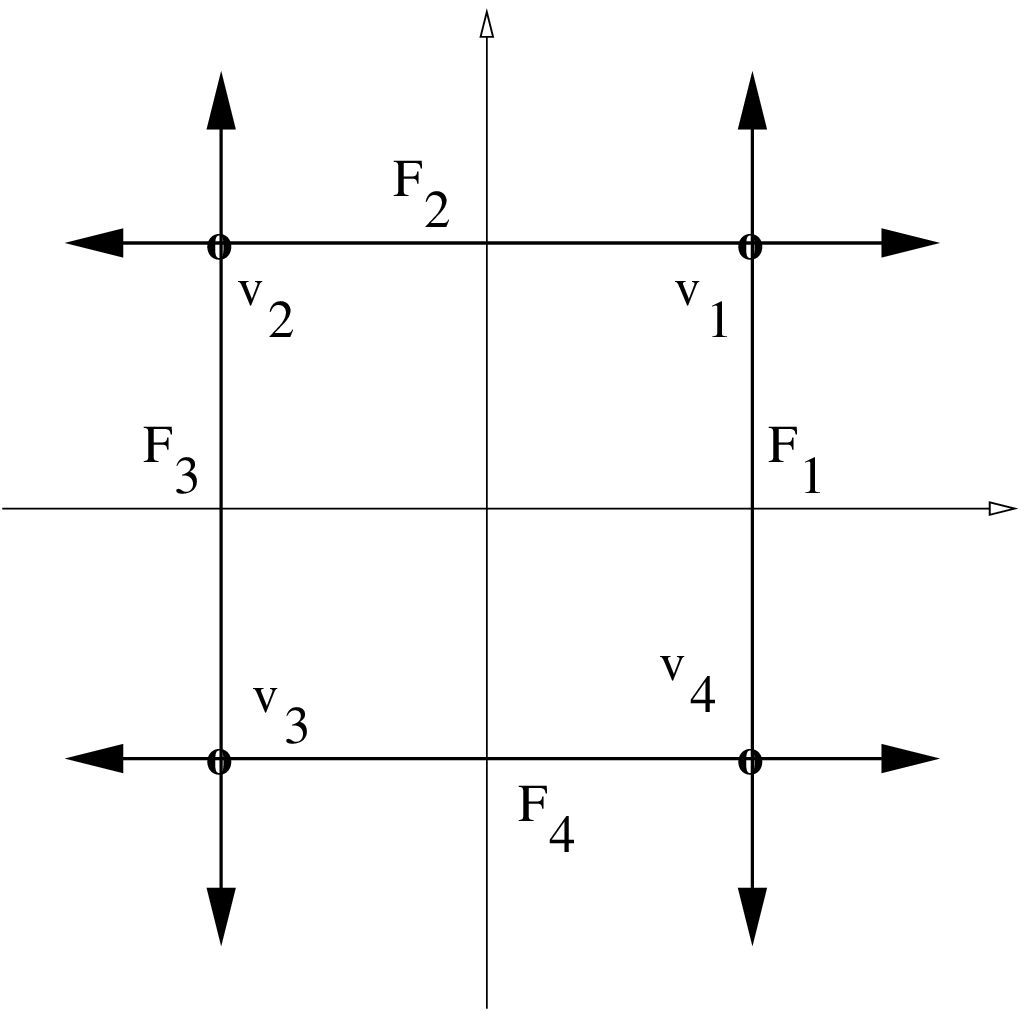}} 
\caption{Polytope normal to the fan of $\c P^1 \times \c P^1$ and the vecotrs $u_{\sigma, \rho}$.}
~\label{polytope2}
\end{figure}

Let $\gamma = (\gamma_1, \gamma_2)$ be a $1$-parameter subgroup.
We get the following table for the values of the $f_\rho$:
$$
\begin{array}{|c|c|c|c|c|}
\hline & & & & \\
& z_1 & z_2 & z_3 & z_4 \\ \hline & & & & \\
f_1 & \gamma_1 & 0 & 0 & \gamma_1 \\ \hline & & & & \\
f_2 & \gamma_2 & \gamma_2 & 0 & 0 \\ \hline & & & & \\
f_3 & 0 & -\gamma_1 & -\gamma_1 & 0 \\ \hline & & & & \\
f_4 & 0 & 0 & -\gamma_2 & -\gamma_2 \\ \hline
\end{array}
$$
and hence, $\Z = \{ (\gamma_1, \gamma_2, 0, 0), (0, \gamma_2,
-\gamma_1, 0), (0, 0, -\gamma_1, -\gamma_2), (\gamma_1, 0, 0,
-\gamma_2) \} \subset \r^4$. Note that the points in $\Z$ lie on
the same $2$-plane defined by $f_1 - f_3 = \gamma_1$ and $f_2 -
f_4 = \gamma_2$. Also, no three of them are colinear. Thus, one
can see that $Gr_0 \cong \c, Gr_1 \cong \c ^2, Gr_2 \cong \c$ and
$Gr_i = \{0\}, i>2$. If $\{x,y\}$ is a basis for $Gr_1 A(\Z)$
then, $H^*(\c P^1 \times \c P^1, \c) \cong Gr A(\Z) \cong \c[x,y]
/ \langle x^2, y^2 \rangle$. In fact, any set of $4$ points lying
on the same $2$-plane such that no three are colinear can give the
cohomology of $\c P^1 \times \c P^1$.
\end{Ex}

\noindent University of British Columbia, Vancouver, B.C. \\{\it
Email address:} {\sf kaveh@math.ubc.ca}

\end{document}